\documentclass[12pt,]{amsart}
\usepackage{amsmath}
\usepackage{amscd}
\usepackage{amssymb}
\usepackage{latexsym}
\input xy
\xyoption{all} \CompileMatrices
\newtheorem{theorem}{Theorem}[section]
\newtheorem{lemma}[theorem]{Lemma}
\newtheorem{corollary}[theorem]{Corollary}
\newtheorem{proposition}[theorem]{Proposition}

\newtheorem{remark}[theorem]{Remark}
\newtheorem{definition}[theorem]{Definition}

\newcommand{\bgl}{\begin{equation}}         
\newcommand{\egl}{\end{equation}}
\newcommand{\bgln}{\begin{eqnarray}}        
\newcommand{\egln}{\end{eqnarray}}
\newcommand{\bglnoz}{\begin{eqnarray*}}     
\newcommand{\eglnoz}{\end{eqnarray*}}
\newcommand{\btheo}{\begin{theorem}}
\newcommand{\etheo}{\end{theorem}}
\newcommand{\blemma}{\begin{lemma}}
\newcommand{\elemma}{\end{lemma}}
\newcommand{\bproof}{\begin{proof}}
\newcommand{\eproof}{\end{proof}}
\newcommand{\bbew}{\begin{beweis}}
\newcommand{\ebew}{\end{beweis}}
\newcommand{\bremark}{\begin{remark}\em}
\newcommand{\eremark}{\end{remark}}
\newcommand{\bdefin}{\begin{definition}}
\newcommand{\edefin}{\end{definition}}
\newcommand{\bprop}{\begin{proposition}}
\newcommand{\eprop}{\end{proposition}}
\newcommand{\bcor}{\begin{corollary}}
\newcommand{\ecor}{\end{corollary}}
\newcommand{\mn}{\par\medskip\noindent}

\newcommand{\cD}{\mathcal D}

\newcommand{\cL}{\mathcal L}

\newcommand{\cO}{\mathcal O}

\newcommand{\cQ}{\mathcal Q}

\newcommand{\lori}{\longrightarrow}

\newcommand{\sgn}{\rm sgn}
\newcommand{\ve}{\varepsilon}
\newcommand{\vp}{\varphi}

\def\SEMI{\mbox{$\times\kern-2pt\vrule height5pt width.6pt \kern3pt $}}

\newcommand{\Hom}{{\rm Hom}\,}

\newcommand{\Ker}{{\rm Ker\,}}
\newcommand{\Img}{{\rm Im\,}}

\newcommand{\Coker}{{\rm Coker}\,}
\newcommand{\id}{{\rm id}}

\newcommand{\orb}{{\rm orb\,}}

\def\Cz{\mathbb{C}}
\def\Fz{\mathbb{F}}

\def\Nz{\mathbb{N}}

\def\Rz{\mathbb{R}}
\def\Zz{\mathbb{Z}}

\def\Tz{\mathbb{T}}
\def\A{\mathfrak{A}}
\def\Af{\mathfrak{A}[\varphi]}
\def\Afp{\mathfrak{A}[\varphi /\psi]}
\newcommand{\abs}[1]{\lvert#1\rvert}     

\begin{document}

\title[Endomorphisms and polymorphisms
of compact abelian groups]{C*-algebras associated with endomorphisms
and polymorphisms of compact abelian groups}
\date{\today}
\author[J. Cuntz]{Joachim Cuntz$^1$}
\author[A. Vershik]{Anatoly Vershik$^2$}
\address{Joachim Cuntz, Mathematisches Institut, Einsteinstr.62, 48149
M\"unster, Germany\\Anatoly Vershik, St.Petersburg Department of
Steklov Institute of Mathematics, 27 Fontanka, St.Petersburg 191023,
Russia}\email{cuntz@uni-muenster.de, vershik@pdmi.ras.ru}
\thanks{{$^1$}Research supported by DFG through
CRC 878 and by ERC through AdG 267079,{$^2$}Research supported by a
Humboldt Research Award and by RFBR grant 11-01-12092-ofi$_m$ }
\begin{abstract}\noindent
A surjective endomorphism or, more generally, a polymorphism in the
sense of \cite{SV}, of a compact abelian group $H$ induces a
transformation of $L^2(H)$. We study the C*-algebra generated by
this operator together with the algebra of continuous functions
$C(H)$ which acts as multiplication operators on $L^2(H)$. Under a
natural condition on the endo- or polymorphism, this algebra is
simple and can be described by generators and relations. In the case
of an endomorphism it is always purely infinite, while for a
polymorphism in the class we consider, it is either purely infinite
or has a unique trace. We prove a formula allowing to determine the
$K$-theory of these algebras and use it to compute the $K$-groups in
a number of interesting examples.
\end{abstract}\maketitle

\section{Introduction}
Let $H$ be a compact abelian group. Let $\alpha$ be an automorphism,
a surjective endomorphism or an algebraic polymorphism (see below)
of $H$. For Haar measure on $H$, the transformation $\alpha$ will
define an operator $s_\alpha$ on the Hilbert space $L^2H$. This
operator will be unitary, isometric or a Markov operator,
respectively. We can form the C*-algebra $C^*(s_\alpha,C(H))$
generated in $\cL(L^2H)$ by $s_\alpha$ together with $C(H)$ acting
as multiplication operators on $L^2H$. The case where $\alpha$ is an
automorphism is classical. In this case $s_\alpha$ is unitary and
$C^*(s_\alpha,C(H))$ is essentially the crossed product by $\alpha$.
The study of the case where $\alpha$ is an endomorphism was started
in the 70s and by now there is quite some literature. The best known
example is the case where $\alpha$ is the shift on
$H=\prod_{k\in\Nz}\Zz/n$. The corresponding C*-algebra is $\cO_n$,
\cite{CuCMP}. This case was also considered in \cite{VerAr} from the
point of view of W*-algebras and factor representations.

The notion of a general measure preserving polymorphism was
suggested in \cite{VerPo1},\cite{VerPo2}. The case of algebraic
polymorphisms was discussed in \cite{SV}. The question of studying
the C*-algebra corresponding to a polymorphism has hardly been
touched upon so far, cf. however \cite{ExR} where C*-algebras
associated with such polymorphisms have been discussed in the
setting of an associated groupoid.

In the present paper we analyze the structure of such C*-algebras
and study their $K$-theory. We start with the case of an algebraic
endomorphism $\alpha$ of the compact abelian group $H$. We assume
that $\alpha$ is surjective with finite kernel and exact (i.e. the
union of the kernels of the $\alpha^n$ is dense, see section
\ref{end}).

In fact, to study $C^*(s_\alpha,C(H))$, most of the time it is more
useful to work, rather than with $\alpha$, with the dual
endomorphism $\vp=\hat{\alpha}$ of the dual group $G=\hat{H}$.
Fourier transform transforms $C^*(s_\alpha,C(H))$ isomorphically
into the C*-algebra $\Af$ acting on $\ell^2G$. As in \cite{Hirsh}
and \cite{CuLi}, but still somewhat surprisingly, this C*-algebra
$\Af$ which is originally defined by a concrete representation, can
also be characterized as a universal algebra given by generators and
relations. The structure of $\Af$ (and thus of $C^*(s_\alpha,C(H))$)
is governed by two ``complementary'' maximal abelian subalgebras,
one being the algebra of continuous functions on $H$, the other one
the algebra of continuous functions on a compactification (with
respect to $\vp$) of $G$. The algebras $\Af$ all have a similar
structure, in particular they are simple, nuclear and purely
infinite. They therefore belong to a very well understood class of
C*-algebras. In particular by the Kirchberg-Phillips classification
\cite{Kip,Phi}, they are completely determined by their $K$-theory.

In section \ref{Kend} we derive a Pimsner-Voiculescu type formula
that can be used to determine the $K$-theory of $\Af$. We prove that
there is an exact sequence of the form

\bgl\label{PV0}\xymatrix{K_*C(H)\ar[r]^{1-b(\vp)}&
K_*C(H)\ar[r]&K_*\Af\ar@/^5mm/[ll]}\egl\mn

We should point out that $b(\vp)$ is not simply the map induced by
$\hat{\vp}$, even though it is related to this map by a simple
equation. The determination of $b(\vp)$ in specific examples
sometimes requires extra work. Since $K_*C(H)$ is always
torsion-free, the exact sequence (\ref{PV0}) is particularly useful
for computations. We use it to explicitly determine the $K$-theory
of $\Af$ for several examples, including the case of endomorphisms
of $\Tz^n$, of $\prod_k\Zz/n$ and of a solenoid group.

Algebras such as those in section \ref{end} have been studied by
quite a few authors. The simplicity of the algebra $\Af$ and its
description as a universal algebra has been established already in
\cite{Hirsh} even in a more general setting. Constructions along the
same lines are considered in the thesis of F.Vieira, \cite{Vie}. As
pointed out to us by R.Exel, the simplicity of $\Af$ could also be
established using an approach as in \cite{Exel} and \cite{ExV}. One
virtue of our approach here is its simplicity together with the fact
that it reveals interesting structural properties of $\Af$ and its
canonical subalgebras.

Special cases of the algebra $\Af$ for $H=\Tz^n$ or for
$H=\prod_{k}\Zz/p$ had also occurred before in
\cite{CuLi,CuLiX,CuLiP}, where again it was shown that in these
examples $\Af$ is purely infinite simple and its $K$-theory was
partially computed. Our proof here that $\Af$ is purely infinite
simple is very similar to that in \cite{CuLi}. For the case of an
endomorphism of $\Tz^n$, an algebra which is easily seen to be
isomorphic to $\Af$ has been described as a Cuntz-Pimsner algebra in
\cite{EHR}, using Exel's concept of a transfer operator \cite{Exel}.
In this paper, it was also proved that, for an expansive
endomorphism of $\Tz^n$, the algebra is simple purely infinite and
its $K$-theory was determined (using Pimsner's extension which leads
to a sequence similar to (\ref{PV0})).  Our computation of the
$K$-theory is somewhat simpler and more general.

Let us now turn to the case of an algebraic polymorphism of the
compact abelian group $H$. This is a multivalued map determined
typically by a pair of endomorphisms. We will restrict ourselves to
the case of what we call rational polymorphisms.

We introduce a concept of independence (motivated by the notion of
relatively prime principal ideals in number theory) for two
commuting surjective endomorphisms of a compact abelian group. This
concept is interesting in itself. It leads to a very special form of
the orbit equivalence relation for the semigroup generated by the
two endomorphisms. We discuss this point briefly in section
\ref{orbit}. The concept of independence could also be generalized
to more general commuting pairs of not necessarily algebraic
endomorphisms of a suitable measure space.

A rational polymorphism is then roughly speaking a product of
$\alpha$ by the ``inverse'' of $\beta$ for two commuting and
independent endomorphisms $\alpha$ and $\beta$ of the compact
abelian group $H$. In \cite{SV} a Markov operator on $L^2(H)$ was
associated with an algebraic polymorphism. This operator fits very
well with the construction mentioned above. For a rational
polymorphism it is a partial isometry. We can then generalize the
construction and analysis of the algebra $\Af$ to a similar
construction of an algebra $\Afp$ associated with a rational
polymorphism induced by a pair of commuting and independent
endomorphisms (actually below $\vp$ and $\psi$ will denote the dual
endomorphisms of the dual group). Again, these algebras can be
characterized by generators and relations and are simple and
nuclear. They are purely infinite if the kernels of the given
endomorphisms of the compact group $H$ do not have the same number
of elements. If these kernels have the same size, then $\Afp$ has a
unique trace.

We also study the $K$-theory of $\Afp$. By an argument similar to
the case of a single endomorphism, we show that the $K$-groups
satisfy an exact sequence of Pimsner-Voiculescu type but with
somewhat more complicated ingredients. In particular the map
$1-b(\vp)$ in formula (\ref{PV0}) is replaced by $b(\psi) -b(\vp)$
where $b(\psi)$ is the map corresponding to the second endomorphism
$\psi$. Again, this sequence often suffices to explicitly compute
the $K$-theory of $\Afp$.

In all our arguments concerning the structure and the $K$-theory of
$\Afp$, the independence of the pair $(\vp,\psi )$ plays a crucial
role.

\section{The algebra associated with an endomorphism of a compact
abelian group}\label{end}

Let $H$ be a compact abelian group and $G=\hat{H}$ its dual discrete
group. We usually denote the group operation on $H$ by
multiplication with neutral element 1 and on $G$ by addition with
neutral element 0. We also assume that $G$ is countable. Let
$\alpha$ be a surjective endomorphism of $H$ with finite kernel. We
denote by $\vp$ the dual endomorphism $\chi\mapsto \chi\circ \vp$ of
$G$ (i.e. $\vp=\hat{\alpha}$). By duality, $\vp$ is injective and
has finite cokernel, i.e. the quotient $G/\vp G$ will be finite.
Both $\alpha$ and $\vp$ induce isometric endomorphisms $s_\alpha$
and $s_\vp$ of the Hilbert spaces $L^2H$ and $\ell^2G$,
respectively.

We will also assume that $$\bigcap_{n\in\Nz}\vp^nG=\{0\}$$ which, by
duality, means that $$\bigcup_{n\in\Nz}\Ker \alpha^n$$ is dense in
$H$ (this implies in particular that $H$, $G$ can not be finite).

\bremark\label{exact} In ergodic theory a measure preserving
endomorphism $T:X\to X$ is called \emph{exact} if it has the
property that $\bigcap_{n\in\Nz}T^{-n}(\mathfrak M)=\mathfrak N$
where $\mathfrak M$ denotes the sigma-algebra of all measurable sets
while $\mathfrak N$ is the trivial sigma-algebra of sets of measure
0 or 1. For an algebraic endomorphism of the compact group $H$ this
condition means exactly that the subgroup $\bigcup_{n\in\Nz}\Ker
\alpha^n$ is dense in $H$. Exact endomorphisms of $H$ can be
characterized as those endomorphisms that have no quotient factor
automorphisms. Among all endomorphisms the exact ones are generic
and are the most interesting ones.\eremark

The standing assumption for the rest of the paper will be that
$\alpha$ is a surjective endomorphism of the compact abelian group
$H$ with non-trivial finite kernel, satisfying the exactness
condition of Remark \ref{exact} (i.e. the union of the kernels of
$\alpha^n$ is dense), or equivalently, that the dual endomorphism
$\vp$ of the dual group $G$ is injective with $1< \abs{G/\vp
(G)}<\infty$ and $\bigcap_{n\in \Nz}\vp^n(G)=\{0\}$.

In this section we are going to describe the C*-algebra
$C^*(s_\alpha,C(K))$ generated in $\cL(L^2H)$ by $C(K)$, acting by
multiplication operators, and by the isometry $s_\alpha$. Via
Fourier transform it is isomorphic to the C*-algebra
$C^*(s_\vp,C^*G)$ generated in $\cL(\ell^2 G)$ by $C^*G$, acting via
the left regular representation, and by the isometry $s_\vp$. These
two unitarily equivalent representations are useful for different
purposes.

\bremark It seems to be interesting to replace in this construction
Haar measure on $H$ by another $\alpha$-invariant measure $\mu$ and
to consider then the C*-algebra of operators on $L^2(H,\mu)$
generated by $s_\alpha$ and $C(K)$. In this way one would obtain for
instance algebras analogous to the algebras $\cO_A$ introduced in
\cite{CuKr}.\eremark

$C^*(s_\vp,C^*G)$ is generated by an isometry $s=s_\vp$ together
with unitary operators $u_g,\, g\in G$, satisfying the relations

\bgl\label{rel}u_gu_h=u_{gh}\qquad s u_g=u_{\vp (g)}s\qquad
\sum_{g\in G/\vp G}u_gss^*u_g^*=1\egl

Even though the analysis of the structure of $C^*(s_\vp,C^*G)$ is
essentially a straightforward generalization of constructions in
\cite{Cun} and \cite{CuLi}, it is interesting enough to merit a
separate discussion.

\bdefin Let $H,G$ and $\alpha, \vp$ be as above. We denote by $\Af$
the universal C*-algebra generated by an isometry $s$ and unitary
operators $u_g,\, g\in G$ satisfying the relations
(\ref{rel}).\edefin

We will show that $\Af\cong C^*(s_\alpha,C(H))\cong
C^*(s_\vp,C^*G)$.

\blemma\label{D} The C*-subalgebra $\cD$ of $\Af$ generated by all
projections of the form $u_gs^ns^{*n}u_g^*$, $g\in G, n\in \Nz$ is
commutative. Its spectrum is the ``$\vp$-adic completion''
$$G_\vp=\mathop{\lim}\limits_{ {\scriptstyle\longleftarrow_n} }G/\vp^nG$$
It is an inverse limit of the finite spaces $G/\vp^nG$ and becomes a
Cantor space with the natural topology.

$G$ acts on $\cD$ via $d\mapsto u_gdu_g^*$, $g\in G,\, d\in \cD$.
This action corresponds to the natural action of the dense subgroup
$G$ on its completion $G_\vp$ via translation. The map $\cD\to\cD$
given by $x\mapsto sxs^*$ corresponds to the map induced by $\vp$ on
$G_\vp$.\elemma

\bproof It is easily checked that the maps
$C^*(\{u_gs^ns^{*n}u_g^*:g\in G\})\to C(G/\vp^n(G))$, $n\in\Nz$ that
map $u_gs^ns^{*n}u_g^*$ to the characteristic function of the one
point set $\{\vp^n(G)+g\}$ extend to an isomorphism of the inductive
limits with the asserted properties (this is basically the same
situation as in \cite{Cun} or in \cite{CuLi}).\eproof

From now on we will denote the compact abelian group $G_\vp$ by $K$. By construction, $G$ is a dense subgroup of $K$.
The dual group of $K$ is the discrete abelian group
$$L = \mathop{\lim}\limits_{ {\scriptstyle\longrightarrow_n}}\Ker
(\alpha^n :H\to H)$$ Because of the condition that we impose on
$\alpha$, $L$ can be considered as a dense subgroup of $H$.

The groups $K$ and $L$ will play an important role in the analysis
of $\Af$. They are in a sense complementary to $H$ and $G$. By Lemma
\ref{D}, the C*-algebra $\cD$ is isomorphic to $C(K)$ and to
$C^*(L)$.

\blemma\label{B} The C*-subalgebra $B_\vp$ of $\Af$ generated by
$C(H)$ together with $C(K)$ (or equivalently by $C^*G$ together with
$C^*L$) is isomorphic to the crossed product $C(K)\rtimes G$. It is
simple and has a unique trace. \elemma \bproof The action of the
dense subgroup $G$ by translation on $K$ is obviously minimal (every
orbit is dense). Therefore the crossed product $C(K)\rtimes G$ is
simple. It also has a unique trace, the Haar measure on $K$ being
the only invariant measure. The fact that an invariant measure on
$K$ extends uniquely to a trace on the crossed product, if all the
stabilizer groups are trivial, is well known, but not easy to pin
down in the literature. Here is a very simple argument in the
present case: Let $E:C(K)\rtimes G\to C(K)$ be the canonical
conditional expectation and let $e_1^{(n)},e_2^{(n)},\ldots
e_{N(\vp^n)}^{(n)}$, with $N(\vp^n)=\abs{G/\vp^n(G)}$, be the
minimal projections in $C(G/\vp^n(G))\subset C(K)$. Then, for any
$x$ in the crossed product,
$\sum_{i=1}^{N(\vp^n)}e_i^{(n)}xe_i^{(n)}$ converges to $E(x)$ for
$n\to\infty$. For any trace $\tau$ on the crossed product, we have
$$\tau (x)=\tau\Big(\sum_{i=1}^{N(\vp^n)}e_i^{(n)}x\Big)=\tau
\Big(\sum_{i=1}^{N(\vp^n)}e_i^{(n)} xe_i^{(n)}\Big)$$ and therefore
$\tau (x)=\tau (E(x))$.

Finally, by Lemma \ref{D}, $B_\vp$ is generated by a covariant
representation of the system $(C(K),\,G)$. The induced surjective
map $C(K)\rtimes G\to B_\vp$ has to be injective, thus an
isomorphism.\eproof

The simple algebra $B_\vp$ has a very interesting structure. In
particular it can be represented as a crossed product in different
ways. To see this, note first that $C(K)\rtimes G$ can also be
written as $C^*(L)\rtimes G$ ($L$ being the dual group of $K$). The
action of $G$ on $C^*(L)$, dual to the action of $G$ on $C(K)$ by
translation, is determined by the commutation relation

\bgl\label{comm}u_gw_l=\langle g|\,l\rangle\, w_lu_g\,\quad g\in
G,\,l\in L\egl

Here we denote the unitary generators of $C^*(G)$ and $C^*(L)$ by
$u_g$, $w_l$, respectively and use the embedding of $L$ into the
dual group $H$ of $G$ to obtain the pairing between $G$ and $L$.
Therefore $B_\vp$ can also be described as the universal C*-algebra
generated by unitary representations $g\mapsto u_g$ and $l\mapsto
w_l$ of $G$ and $L$ satisfying the commutation relation
(\ref{comm}). In particular, the relation (\ref{comm}) can also be
interpreted as an action of $L$ on $C^*(G)$ and we see that $B_\vp$
is also isomorphic to the crossed product $C^*(G)\rtimes L$.

Now similarly the action of $L$ on $C^*(G)$ determined by the
commutation relation (\ref{comm}) corresponds to the action of $L$
on $C(H)$ by translation, under duality. Therefore we get the
following intriguing chain of isomorphisms

$$B_\vp\cong C(K)\rtimes G\cong C^*L\rtimes G\cong C^*G\rtimes L\cong
C(H)\rtimes L$$

This chain of isomorphisms is similar - though easier - to the
duality result in \cite{CuLiX}.\mn

The map $x\mapsto sxs^*$ defines a natural endomorphism $\gamma_\vp$
of $B_\vp$.

\btheo\label{pure} The algebra $\Af$ is simple, nuclear and purely
infinite.

Moreover, it is isomorphic to the semigroup crossed product
$B_\vp\rtimes_{\gamma_\vp} \Nz$ (i.e. to the universal unital
C*-algebra generated by $B_\vp$ together with an isometry $t$ such
that $txt^*=\gamma_\vp (x)$, $x\in B_\vp$).\etheo

\bproof $\Af$ contains $B_\vp$ as a unital subalgebra. The condition
$\alpha_\lambda (s)=\lambda s,\,\alpha_\lambda (b)= b,\, b\in B_\vp$
defines for each $\lambda\in\Tz$ an automorphism of $A$ and
integration of $\alpha_\lambda (x)$ over $\Tz$ determines a faithful
conditional expectation $A\to B_\vp$. The proof now is very similar
to the corresponding proof in \cite{CuLi}. The representation as a
crossed product $C(K)\rtimes G$ of $B_\vp$ gives a natural faithful
conditional expectation $B_\vp\to \cD\cong C(K)$. The composition of
these expectations gives a faithful conditional expectation
$E:\Af\to \cD$.

Now, this expectation can be represented in a different way using
only the internal structure of $\Af$. The relations (\ref{rel})
immediately show that the linear combinations of elements of the
form $z=s^{*n}du_gs^m$, $n,m\in\Nz,g\in G,d\in\cD$ are dense in
$\Af$. For such an element $z$ we have $E(z)=s^{*n}ds^n$ if
$n=m,g=0$, and $E(z)=0$ otherwise.

The subalgebra $\cD$ of $\Af$ is the inductive limit of the
finite-dimensional subalgebras $\cD_n\cong C(G/\vp ^nG)$. Note that
the minimal projections in $\cD_n$ are all of the form
$u_gs^ns^{*n}u_g^*$. Let

$$z = d + \sum_{i=1}^m s^{*k_i}d_iu_{g_i}s^{l_i}$$

be an element of $\Af$ such that for each $i$, $k_i\neq l_i$ or
$g_i\neq e$ and such that $d,d_i\in \cD_n$ for some large $n$ (such
elements are dense in $\Af$).

Let also $n$ be large enough so that the projections $u_{g_i}e
u_{g_i}^*$, $i=1,\ldots ,m$, are pairwise orthogonal for each
minimal projection $e$ in $\cD_n$ (this means that the $g_i$
are pairwise distinct mod $\vp^n G$).

We have $E(z)=d$ and there is a minimal projection $e$ in $\cD_n$
such that $E(z)e=\lambda e$ with $\abs{\lambda}=\|E(z)\|$. Since $E$
is faithful, $\lambda>0$ if $z$ is positive $\neq 0$.

Let $h\in G$ such that $e =u_hs^ns^{*n}u_h^*$. Then the product
$e s^{*k_i}d_iu_{g_i}s^{l_i}e$ is non-zero only if
$g_i\vp^{l_i}(h)= \vp^{k_i}(h)$ or $g_i=
\vp^{k_i}(h)\vp^{l_i}(h)^{-1}$ mod $\vp^n(G)$.

Let $f\in \vp^nG$ such that $\vp^{k_i}(f)\neq \vp^{l_i}(f)$ for all
$i$ for which $k_i\neq l_i$ (such an $f$ obviously exists) and let
$k\geq 0$ such that $\vp^{k_i}(f)\neq \vp^{l_i}(f)\mod \vp^{n+k}G$
for those $i$.

Then, setting $h'=hf$, we obtain

$$g_i\vp^{l_i}(h')\neq \vp^{k_i}(h')\mod
\vp^{n+k}(G)\, ,\;i=1,\ldots ,m$$

If we now set $e'=u_{h'}s^{n+k}s^{*(n+k)}u_{h'}^*$, then $e'$ is a
minimal projection in $\cD_{n+k}$, $e'\leq e$ and $e'
s^{*k_i}d_iu_{g_i}s^{l_i}e'=0$ for $i=1,\ldots ,m$.

Every positive element $x\neq 0$ of $\Af$ can be approximated up to
an arbitrary $\ve$ by a positive element $z$ as above. Thus, if
$\ve$ is small enough, $e'xe'$ is close to $\lambda e'$ and
therefore invertible in $e'\Af e'$. Thus the product
$s^{*n+k}u_{h'}^*xu_{h'}s^{n+k}$ is invertible in $\Af$. This shows,
at the same time, that $\Af$ is purely infinite and simple.
Moreover, it follows that the natural map from $\Af$ to the
semigroup crossed product $B_\vp\rtimes_{\gamma_\vp}\Nz$, is an
isomorphism. The fact that this crossed product is nuclear ($B_\vp$
is nuclear and, using a standard dilation,
$B_\vp\rtimes_{\gamma_\vp}\Nz$ is Morita equivalent to a crossed
product $B^\infty_\vp\rtimes_{\gamma^\infty_\vp}\Zz$, where
$B^\infty_\vp$ is nuclear) then shows that $\Af$ is nuclear. \eproof

\bcor The natural map induces an isomorphism $\Af\cong
C^*(s_\vp,C^*G)\cong C^*(s_\alpha,C(K))$.\ecor

\bproof Since $C^*(s_\vp,C^*G)$ is generated by elements satisfying
the relations (\ref{rel}), there is a natural surjective map $\Af\to
C^*(s_\vp,C^*G)\cong C^*(s_\alpha,C(K))$. By simplicity of $\Af$
this map has to be injective.\eproof

\subsection{Examples} Let us now look at a number of examples.

\subsubsection{}\label{On} Let $H=\prod_{k\in \Nz}\Zz/n$,
$G=\bigoplus_{k\in \Nz}\Zz/n$ and $\alpha$ the one-sided shift on
$H$ defined by $\alpha ((a_k))=(a_{k+1})$.

We obtain $K=\prod_{k\in \Nz}\Zz/n=H$ and $L=\bigoplus_{k\in
\Nz}\Zz/n\cong G$. The algebra $B_\vp$ is a UHF-algebra of type
$n^\infty$ and $\Af$ is isomorphic to $\cO_n$. It is interesting to
note that the UHF-algebra $B_\vp$ is generated by two maximal
abelian subalgebras both isomorphic to $C(K)$.

\subsubsection{} Let $H=\Tz$, $G=\Zz$ and $\alpha$ the endomorphism of
$H$ defined by $\alpha (z)=z^n$. The algebra $B_\vp$ is a
Bunce-Deddens-algebra of type $n^\infty$ and $\Af$ is isomorphic to
a natural subalgebra of the algebra $\cQ_\Nz$ considered in
\cite{Cun}. In this case, we also get for $B_\vp$ the interesting
isomorphism $C(\Zz_n)\rtimes \Zz\cong C(\Tz)\rtimes L$ where $\Zz$
acts on the Cantor space $\Zz_n$ by the odometer action (addition of
1) and $L$ denotes the subgroup of $\Tz$ given by all $n^k$-th roots
of unity, acting on $\Tz$ by translation.

\subsubsection{} Let $H=\Tz^n$, $G=\Zz^n$ and $\alpha$ an endomorphism
of $H$ determined by an integral matrix $T$ with non-zero
determinant. We assume that the condition
$$\bigcap_{n\in\Nz}\vp^nG=\{0\}$$ is satisfied (this is in fact not very restrictive).

The algebra $B_\vp$ is a higher-dimensional analogue of a
Bunce-Deddens-algebra. In the case where $H$ is the additive group
of the ring $R$ of algebraic integers in a number field of degree
$n$ and the matrix $T$ corresponds to an element of $R$, the algebra
$\Af$ is isomorphic to a natural subalgebra of the algebra $\A [R]$
considered in \cite{CuLi}. It is also isomorphic to the algebra
studied in \cite{EHR}.

\subsubsection{} Let $H=\Fz_p[[t]]$, $G=\Fz_p[t]$ and $\vp$ an
endomorphism of $G$ determined by multiplication by a non-zero
element $P$ in the ring $\Fz_p[t]$. In the simplest case where
$P=t$, we are back in the situation of \ref{On}. Then the algebra
$B_\vp$ is a UHF-algebra of type $p^\infty$ and $\Af$ is isomorphic
to $\cO_p$.

For endomorphisms $\alpha,\vp$ induced by an arbitrary $P$, the
algebra $\Af$ is naturally a subalgebra of the ring C*-algebra
considered in \cite{CuLiP}.

More generally, instead of multiplication by $P$, we could also
consider an arbitrary $\Fz_p$-linear injective endomorphism $\vp$ of $\Fz_p[t]$ with
finite-dimensional cokernel and such that the intersection of the images of all powers of $\vp$ is 0.

\subsubsection{}\label{sol} Let  $p$ and $q$ be natural numbers that
are relatively prime and $\gamma$ the endomorphism of $\Tz$ defined
by $z\mapsto z^p$. We take

$$H=\mathop{\lim}\limits_{ \mathop{{\scriptstyle\longleftarrow}}
\limits_\gamma}\Tz\quad\qquad G=\Zz[\frac {1}{p}]$$\mn

$\alpha_q$ the endomorphism of $H$ induced by $z\mapsto z^q$ and
$\vp_q$ the endomorphism of $G$ defined by $\vp_q(x)=qx$. These
endomorphisms satisfy our hypotheses ($G/\vp^nG=\Zz/q^n$ and
$\bigcap \vp^nG=\{0\}$). We find that $K=\Zz_q$ (the q-adic
completion of $\Zz$).

\section{Computation of the $K$-theory for $\Af$}\label{Kend}

In section \ref{end} we had considered the natural subalgebra
$B_\vp=C(K)\rtimes G$ of $\Af$. Since, by definition,
$$K=\mathop{\lim}\limits_{ {\scriptstyle\longleftarrow_n}}G/\vp^nG$$
we can represent $B_\vp$ as an inductive limit
$B=\mathop{\lim}\limits_{ {\scriptstyle\longrightarrow_n}}B_n$ with
$B_n=C(G/\vp^nG )\rtimes G$.

It is well known (``imprimitivity'') that $$C(G/\vp^nG )\rtimes
G\cong M_{N(\vp)}(C^*(\vp^nG))$$ where $N(\vp)=|G/\vp^nG|$. Consider
the natural inclusion

$$C^*G\cong C^*(\vp^n G)\lori M_{N(\vp)}(C^*(\vp^n G))\cong B_n$$

into the upper left corner of $M_{N(\vp )}$, considered as a map
$C^*G\to B_n$. This map induces an isomorphism $\kappa_n
:K_*(C^*G)\to K_*(B_n)$ in $K$-theory.

We let moreover $\iota_n$ denote the map $K_*(B_n)\to K_*(B_{n+1})$
induced by the inclusion $B_n\hookrightarrow B_{n+1}$ and define

$$b(\vp)_n :K_*(C^*(G))\lori K_*(C^*(G))$$

by $b(\vp)_n =\kappa_{n+1}^{-1}\iota_n\kappa_n$.

Now, the commutative diagram

$$\xymatrix{C^*(G)\ar[r]^{\kappa_0}\ar@/_1pc/[rdd]_{\kappa_n}&B_0\ar[d]
^\cong\ar@^{(->}[r]&
B_1\ar[d]_\cong\quad&C^*(G)\ar[l]_{\kappa_1}\ar@/^1pc/[ldd]^{\kappa_{n+1}}
\\
&C^*(\vp^nG)\ar@^{(->}[d]\ar@^{(->}[r]&\ar@_{(->}[d] C(X_n)\rtimes
\vp^n
G\quad&\qquad\quad\\
&B_n\ar@^{(->}[r]& B_{n+1}&}$$

with $X_n=\vp^nG/\vp^{n+1}G$, shows that $b(\vp)_n=b(\vp)_0$ for all
$n$. We write $b(\vp)$ for this common map.

We obtain the following commutative diagram

$$\xymatrix{K_*(C^*(G))\ar[r]^{b(\vp)}\ar@{=}[d]&
K_*(C^*(G))\ar[r]^{b(\vp)}
\ar[d]^{\kappa_1}&K_*(C^*(G))\ar[r]\ar[d]^{\kappa_2}&\\
K_*(B_0)\ar[r]^{\iota_0}&K_*(B_1)\ar[r]^{\iota_1}&K_*(B_2)\ar[r]&}$$

One immediate consequence is the following formula for the
$K$-theory of $B_\vp$

\bgl\label{KBf}K_*(B_\vp)=\mathop{\lim}\limits_{
\mathop{{\scriptstyle\longrightarrow}}
\limits_{b(\vp)}}K_*(C^*(G))\egl

We note however, that the problem remains to determine a suitable
formula for the map $b(\vp)$, given a specific endomorphism $\vp$.
We will consider a few examples below.

Since, by Theorem \ref{pure}, $\Af$ can be represented as a crossed
product $B_\vp\rtimes_{\gamma_\vp}\Nz$, we are now in a position to
derive a formula for the $K$-theory of $\Af$.

\btheo\label{KA} The $K$-groups of $\Af$ fit into an exact sequence
as follows

$$\xymatrix{&K_0C^*(G)\ar[r]^{1-b(\vp)}&K_0C^*(G)\ar[rd]&\\
K_1\Af\ar[ru]& & &
K_0\Af\ar[ld]\\&K_1C^*(G)\ar[lu]&K_1C^*(G)\ar[l]_{1-b(\vp)}& }$$
\etheo

\bproof From Theorem \ref{pure} we know that $\Af$ is isomorphic to
the semigroup crossed product $B_\vp\rtimes_{\gamma_\vp}\Nz$. Using
the Pimsner-Voiculescu sequence \cite{PV} in combination with a
simple dilation argument as in \cite{CuII} (or directly appealing to
the results in \cite{Pim} or in \cite{KS}) we see that there is an
exact sequence

\bgl\label{PV}\xymatrix{K_*B_\vp\ar[r]^{1-\gamma_{\vp*}}&
K_*B_\vp\ar[r]&K_*\Af\ar@/^5mm/[ll]}\egl\mn

In order to determine the kernel and cokernel of the map
$K_*B_\vp\mathop{\lori}\limits^{1-\gamma_{\vp*}} K_*B_\vp$, consider
the commutative diagram

$$\xymatrix{K_*C^*(G)\ar[r]^{b(\vp)}\ar[d]^{\kappa_0}_=&K_*C^*(G)\ar[r]^
{b(\vp)}
\ar[d]^{\kappa_1}&K_*C^*(G)\ar[r]^{\quad b(\vp)}\ar[d]^{\kappa_2} &\\
K_*B_0\ar[r]^{\iota_0}\ar[d]&K_*B_1\ar[r]^{\iota_1}\ar[d]&K_*B_2
\ar[r]^{\iota_2}\ar[d]&\\
K_*B_\vp\ar[r]^{=}&K_*B_\vp\ar[r]^{=}&K_*B_\vp\ar[r]^{=}&}$$

By construction, it is clear that
$\gamma_{\vp*}\kappa_n=\kappa_{n+1}$ (where we still denote the
composition $K_*C^*(G)\mathop{\lori}\limits^{\kappa_n} K_*B_n\to
K_*B_\vp$ by $\kappa_n$). Let $\kappa$ denote the map (isomorphism)
$$\kappa: \mathop{\lim}\limits_{
{\scriptstyle\longrightarrow_{b(\vp)}}}K_*C^*(G)\lori K_*(B_\vp)$$
induced by the commutative diagram.

For an element of the form $[x_0,x_1,\ldots]$ in the inductive limit
we then obtain

$$\gamma_{\vp*}\circ\kappa ([x_0,x_1,x_2,\ldots])= \kappa ([a,x_0,x_1,\ldots])$$
(where $a$ is arbitrary).

Therefore the exact sequence (\ref{PV}) becomes isomorphic to

\bgl\label{di}\xymatrix{\mathop{\lim}\limits_{
{\scriptstyle\longrightarrow_{b(\vp)}}}K_*C^*(G)\ar[r]^{1-\sigma}&
\mathop{\lim}\limits_{
{\scriptstyle\longrightarrow_{b(\vp)}}}K_*C^*(G)\ar[r]&K_*\Af\ar@/^5mm/[ll]}\egl\mn

where $\sigma$ is the shift defined by

$$\sigma ([x_0,x_1,x_2,\ldots])=  [a,x_0,x_1,\ldots]$$

Consider the natural map $j:K_*C^*G\to \mathop{\lim}\limits_{
{\scriptstyle\longrightarrow_{b(\vp)}}}K_*C^*(G)$ defined by
$$j(x)=[x,b(\vp) (x),b(\vp)^2(x),\ldots]$$

If $(1-\sigma)[x_0,x_1,\ldots]=0$, then there is $n$ such that
$x_n=x_{n+1}=b(\vp) (x_n)$ and thus
$[x_0,x_1,\ldots]=[x_n,x_n,\ldots ]$. This shows that $\Ker
(1-\sigma)=j(\Ker (1-b(\vp) ))\cong \Ker (1-b(\vp) )$.

If we divide $\mathop{\lim}\limits_{
{\scriptstyle\longrightarrow_{b(\vp)}}}K_*C^*(G)$ by $\Img
(1-\sigma)$, then $[x_0,x_1,\ldots]$ becomes identified with
$[x_1,x_2,\ldots]$ and thus to an element of the form $[x,b(\vp)
(x),\ldots]$ which is in the image of $j$. Also $j$ maps $\Ker
(1-b(\vp))$ to $\Ker (1-\sigma)$ and thus induces an isomorphism
from the cokernel of $1-b(\vp)$ to the cokernel of $1-\sigma$. This
shows that $j$ induces a transformation from the sequence

$$\xymatrix{K_*C^*(G)\ar[r]^{1-b(\vp)}&
K_*C^*(G)\ar[r]&K_*\Af\ar@/^5mm/[ll]}$$\mn

into the exact sequence (\ref{di}), which is an isomorphism on
kernels and cokernels (in fact $j$ transforms $1-b(\vp)$ not into
$1-\sigma$ but into $1-\sigma^{-1}$ - this however does not affect
exactness). \eproof

There is a connection between $b(\vp)$ and $\vp_*$ which is
described in the following Lemma.

\blemma\label{al} On $K_*C^*(G)$ we have the identity
$b(\vp)\vp_*=N(\vp)\id$ where $N(\vp)=|G/\vp (G)|$.\elemma

\bproof Under the identification $B_1\cong M_{N(\vp)}(C^*(G))$, the
map $\iota_0\kappa_0\vp_*$ is induced by the embedding of
$C^*(G)\cong C^*(\vp G)$ along the diagonal of $M_{N(\vp)}(C^*(G))$.
Therefore $\iota_0\kappa_0\vp_* = N(\vp)\kappa_1$ ($\kappa_1$ is
induced by the embedding in the upper left corner). The assertion
now follows from the definition of $b(\vp)$ as
$\kappa_1^{-1}\iota_0\kappa_0$.\eproof

\subsection{Examples}
\subsubsection{} Let $H=\prod_{k\in \Nz}\Zz/n$, $G=\bigoplus_{k\in
\Nz}\Zz/n$ and $\alpha$ the one-sided shift on $H$ defined by
$\alpha ((a_k))=(a_{k+1})$. As we explained in \ref{On} in this case
$\Af\cong\cO_n$.

It is well known (and easy to see) that $K_0 C(H)\cong
C(H,\Zz),\;K_1C(H)=0$. If we describe the elements of $H$ by
sequences $(x_0,x_1,\ldots)$ with $x_i\in \Zz/n$, then on $f\in
 C(H,\Zz)\cong K_0 C(H)$, the map $b(\vp)$ is given by

$$b(\vp) (f)(x_0,x_1,\ldots)= \sum_{k=0}^{n-1}f(k,x_1,x_2,\ldots )$$

while $\vp_*$ is described by $\vp_*f (x_0,x_1,\ldots)
=f(x_1,x_2,\ldots )$. The application of Theorem \ref{KA} and the
inductive limit description of $K_*B_\vp$ of course leads to the
well known formulas for the $K$-theory of $B_\vp$ and $\cO_n$, i.e
$K_0(B_\vp)=\Zz[\frac {1}{n}],\, K_1(B_\vp)=0$ and
$K_0(\cO_n)=\Zz/(n-1),\, K_1(\cO_n)=0$.

\subsubsection{}

Let $\alpha$ be an endomorphism of $H=\Tz^n$ with finite kernel and
$\vp$ the dual endomorphism of $G=\Zz^n$. We assume that the
intersection of all $\vp^k(\Zz^n)$ is $\{0\}$.

We know that there is a grading (and exterior product) preserving
isomorphism of $K_*(C(\Tz^n))$ with the exterior algebra
$\Lambda^*\Zz^n =\bigoplus_{p=0}^n \Lambda^p\Zz$. The endomorphism
$\vp_*$ of $K_*(C(\Tz^n))$ induced by $\vp$ corresponds to the
endomorphism $\Lambda\vp$ of $\Lambda^*\Zz^n$.

The associated endomorphism $b(\vp)$ of $\Lambda^*\Zz^n$ is
determined by the formula $b(\vp)\vp_*=N(\vp)\,\id$ from Lemma
\ref{al}. In the present case we have $N(\vp )=\abs{\det \vp}$.

Consider the Poincar\'{e} isomorphism $D:\Lambda G\cong \Lambda G'$
(here we write $G'$ for the algebraic dual $\Hom (G,\Zz)$) and
denote by $\vp'$ the endomorphism of $G'$ which is dual to $\vp$
under the natural pairing $G\times G'\to \Zz$).

By \cite{Gr}, (6.62), one has $$\Lambda\vp\, (D\Lambda\vp'
D^{-1})=\det \vp\,\id$$ Therefore the unique solution $b(\vp)$ (in
endomorphisms of $\Lambda\Zz^n$) for the equation
$b(\vp)\,\Lambda\vp=\abs{\det \vp}\;\id$ corresponds under the
Poincar\'{e} isomorphism to $\ve\,\Lambda \vp'$ with $\ve=\sgn (\det
\vp)$. The restriction of $b(\vp)$ to $\Lambda^1\Zz^n\cong\Zz^n$ for
instance is the complementary matrix to $\vp$ determined by Cramer's
rule. Thus we obtain

$$K_*\Af\cong \Lambda G'/(1-\ve\,\Lambda\vp')
\Lambda G'\;\oplus \;\Ker (1-\ve\,\Lambda\vp')$$

where the first term has the natural even/odd grading. The second
term $\Ker (1-\ve\Lambda\vp')$ is $\Lambda^n\Zz^n\cong \Zz$ if $\det
\vp>0$ and $\{0\}$ if $\det\vp<0$. It contributes to $K_0$ if $n$ is
odd and to $K_1$ if $n$ is even.

For instance, if $\vp =k\,\id$, then $\det \vp=k^n$, and $\Lambda
\vp =k^p\,\id$ on $\Lambda^p\Zz^n$. Therefore
$b(\vp)|_{\Lambda^p\Zz^n}=k^{n-p}\,\id$. We thus obtain
$$K_*\Af = \bigoplus_{0\leq p\leq
n}\Lambda^p\Zz^n/(1-k^{n-p})\Lambda^p\Zz^n\;\oplus \;\Zz$$
(understood with the natural even/odd grading for the first term on
the right hand side and $\Zz$ contributing to $K_0$ or $K_1$
depending on the parity of $n$). The same type of formula has been
obtained in \cite{EHR} using an exact sequence similar to the one in
Theorem \ref{KA}, which itself however has been obtained in a rather
different way.

\subsubsection{} Consider the solenoid group $H$ of Example
\ref{sol}, i.e.
$$H=\mathop{\lim}\limits_{ \mathop{{\scriptstyle\longleftarrow}}
\limits_p}\Tz\quad\qquad G=\Zz[\frac {1}{p}]$$\mn with the
endomorphism $\vp_q$ determined on $G$ by $\vp_q(x)=qx$ ($q$ prime
to $p$). The description of $G$ as an inductive limit of groups of
the form $\Zz$ immediately leads to the formulas
$$K_0(C^*G) =\Zz\qquad K_1(C^*G)=\Zz[\frac {1}{p}]$$
Now $\vp_q$ acts as id on $K_0(C^*G)$ and by multiplication by $q$
on $K_1(C^*G)$. Since $N(\vp_q)=q$, we infer from Lemma \ref{al}
that $b(\vp)=q\,\id$ on $K_0(C^*G)$ and $b(\vp)=\id$ on $K_1(C^*G)$.
Thus the exact sequence of Theorem \ref{KA} shows that
$$K_0(\Af)=\Zz/(q-1)+\Zz[\frac {1}{p}]\qquad K_1(\Af)=\Zz[\frac {1}{p}]$$

\section{Algebraic polymorphisms}\label{poly0}

Let $(X,\mu)$ be a Lebesgue space with continuous measure $\mu$ (a
measure space isomorphic to the unit interval with Lebesgue
measure).

\begin{definition}A measure preserving polymorphism $\Pi$ of the
Lebesgue space $(X,\mu)$ to
itself is a diagram between three Lebesgue spaces:
 $$\Pi:\quad (X,\mu)  \longleftarrow (X \times X, \nu)
    \longrightarrow (X, \mu)$$
where the left and right arrows are the projections $\pi_1,\pi_2$
onto the first and second component of the product space $(X \times
X, \nu)$, and where we assume that $\pi_{i*}\nu=\mu$,
$i=1,2$.\end{definition}

Because of the condition $\pi_{i*}\nu=\mu$, the operators
$v_i:L^2(X,\mu)\to L^2(X\times X,\nu),\, i=1,2$, induced by $\pi_i$,
are isometries.

\begin{definition}\label{mark}
The Markov operator $s_\Pi$ corresponding to the polymorphism $\Pi$
is the operator on $L^2(X)$ defined by the formula:  $$s_{\Pi}=v^*_1
v_2,$$.
\end{definition}
The operator $s_\Pi$ is a Markov operator in the sense that it is a
contraction ($\|M\|\leq 1$), it is positive ($Mf\geq 0$ if $f\geq
0$) and it preserves the constant function $1$ ($M(1)=M^*(1)=1$).

\bremark\label{part} A product of the form $v_1^*v_2$ where
$v_1,v_2$ are isometries is a partial isometry, if and only if the
range projection of $v_2$ commutes with the range projection of
$v_1$. Therefore $s_\Pi$ is a partial isometry only under additional
assumptions.\eremark

In this article we will consider only algebraic polymorphisms of
compact groups as discussed in \cite{SV}.

\bdefin[cf. \cite{SV},1.1]  Let $H$ be a compact group which is not
finite. A closed subgroup $C\subset H\times H$ is an algebraic
correspondence of $H$ if $\pi_1 (C)=\pi_2 (C)=H$ for the two
coordinate projections $\pi_i :H\times H\to H$, $i=1,2$.\edefin

Every algebraic correspondence gives rise to a polymorphism
$$\Pi_C:\quad (H,\lambda)  \longleftarrow (H \times H, \nu)
    \longrightarrow (H, \lambda)$$
where $\nu$ is the extension of the Haar measure on $C$ to $H\times
H$ such that $\nu ((H\times H)\backslash C)=0$. Clearly, the
$\pi_i:C\to H$ are homomorphisms and send the Haar measure on $C$ to
the Haar measure on $H$. The associated polymorphism therefore is
measure preserving.

\bdefin\label{alge} Given two surjective endomorphisms $\alpha$ and
$\beta$ of $H$ with finite kernel we define a correspondence
$C_{\alpha,\beta}$ and an associated polymorphism
$\Pi_{C_{\alpha,\beta}}$ by setting
$$C_{\alpha,\beta}= \{(\alpha (h),\beta (h)): h\in H\}$$\edefin

If $\Ker \alpha\cap\Ker\beta=\{1\}$, we can identify $H$ with $C$
via the map $h\mapsto (\alpha (h),\beta (h))$. Under this
identification, the isometries $v_1,v_2$ correspond to the isometric
operators $s_\alpha$ and $s_\beta$ on $L^2(H)$ as in section
\ref{end}. Their range projections commute if $H$ is abelian (under
Fourier transform they simply correspond to the orthogonal
projections onto $\ell^2 (\Img\hat{\alpha})$ and $\ell^2
(\Img\hat{\beta})$). Therefore the associated Markov operator
$s_\Pi$ is a partial isometry if $H$ is abelian and $\Ker
\alpha\cap\Ker\beta=\{1\}$.

\section{Independence}

Let $\varphi$ and $\psi$ be two injective endomorphisms of a
discrete abelian group $G$ satisfying  the conditions of section
\ref{end}.

We will assume that $\vp$ and $\psi$ commute. For much of our
discussion (in particular concerning polymorphisms defined by a pair
of endomorphisms) we need a stronger condition which is described in
the following lemma. This condition is well known in number theory
for the case where $\vp,\psi$ are given by multiplication by an
algebraic number on the ring of algebraic integers.

\blemma\label{prime} Let $\vp,\psi$ be injective commuting
endomorphisms of the abelian group $G$ such that $G/\vp (G)$ and
$G/\psi (G)$ are both finite. Then the following conditions are
equivalent:

\begin{itemize}
  \item [(a)] $\vp (G)+ \psi (G) =G$
  \item [(b)] $\vp (G)\hookrightarrow G$ induces an isomorphism
  $\vp(G)/(\vp(G)\cap\psi(G))\cong G/\psi (G)$.
  \item [(c)] $\vp (G)\cap\psi (G)= \vp\psi (G)$
\end{itemize}\elemma

\bproof The induced map in (b) is injective. It is also surjective
if we assume (a). Thus (a) $\Rightarrow$ (b). The fact that $\vp\psi
(G)\subset \vp (G)\cap\psi (G)$, together with the isomorphism
$G/\psi (G)\cong \vp(G)/\vp\psi (G)$  shows that (b) $\Rightarrow$
(c).

Finally, (c) implies that $(\vp (G)+\psi (G))/(\vp (G)\cap\psi
(G))\cong G/\vp (G)\oplus G/\psi (G)$ $\cong G/(\vp (G)\cap\psi
(G))$ (using the fact that $G/\vp (G)\cong\psi (G)/\psi\vp (G)$).
This implies (a), since both sides of the equation are
finite.\eproof

\bdefin\label{indt}We say that $\vp,\psi$ are independent (or
relatively prime) if the equivalent conditions in \ref{prime} are
satisfied. \edefin Applying condition (a) in \ref{prime}
inductively, we see that then also each power of $\vp$ is prime to
each power of $\psi$ (Proof: $\vp (G)+\psi (G)=\vp (\vp (G)+\psi
(G))+\psi (G)=\vp^2 (G)+\vp\psi(G)+\psi (G) =\vp^2 (G)+\psi (G)$).

\bremark\label{chin} For $\vp,\psi$ independent we have the
following version of the Chinese remainder theorem
$$G/\psi\vp (G) = G/(\psi (G)\cap\vp (G))\cong G/\psi (G)\oplus
G/\vp (G)$$ The second isomorphism is a consequence of condition (a)
above.\eremark

The following lemma is just a reformulation of Lemma \ref{prime} for
the dual group.
\begin{lemma}\label{strong}
Let $\alpha$ and $\beta$ be two commuting surjective endomorphisms
of the compact abelian group $H$ with finite kernel and let
$\vp=\hat{\alpha}$, $\psi=\hat{\beta}$ be the dual endomorphisms of
the dual group $G=\hat{H}$. The following are equivalent
\begin{enumerate}
  \item [(a)] $\vp$ and $\psi$ are independent.
  \item [(b)] $\Ker \alpha \cap\Ker \beta =\{1\}$.
  \item [(c)] $\alpha \,( \Ker \beta) =\Ker \beta$.
  \item [(c')] $\beta\, ( \Ker \alpha) =\Ker \alpha$.
  \item [(d)] The subgroup $\Ker\alpha\Ker\beta$ generated by
  $\Ker\alpha$ and $\Ker\beta$ equals $\Ker\alpha\beta$.
\end{enumerate}
\end{lemma}
\bdefin\label{inp} We say that $\alpha$ and $\beta$ are independent
if they satisfy the equivalent conditions in Lemma \ref{strong}.
\edefin

\bremark\label{mg} Condition (b) in Lemma \ref{strong} says that the
partitions of $H$ into cosets
  with respect to $\Ker\vp$ and $\Ker\psi$ are independent for the Haar
  measure. This means that two functions invariant under $\Ker\vp$ and
  $\Ker\psi$, respectively, are independent in $L^2(H)$ in the
  probabilistic sense. Based on this observation one could also
  define independence for more general (non-algebraic) commuting
  endomorphisms of a measure space.
\eremark

\blemma\label{vw} Let $v$ and $w$ be two isometries  in a unital
C*-algebra. Then the identity $vw^*=w^*v$ implies that $v$ and $w$
commute (but not conversely).\elemma

\bproof The identity $vw^*=w^*v$ implies that $v^*w^*vw =1$. On the
other hand, if a,b are isometries (here b=vw, a=wv) such that
$a^*b=1$, then $a=b$.\eproof

\blemma\label{pi} Let $\vp$ and $\psi$ be injective endomorphisms of
the abelian group $G$. Let $s_\vp$ and $s_\psi$ be the isometries in
  $\cL(\ell^2 G)$ defined by $s_\varphi
(\xi_g)=\xi_{\varphi(g)}$ and $s_\psi (\xi_g)=\xi_{\psi(g)}$. Let
further
  $e_\vp=s_\vp s^*_\vp$ and $e_\psi =s_\psi s^*_\psi$ be the range
  projections. Then
  \begin{itemize}
    \item [(a)] Assume that $\vp$ and $\psi$ commute and are
    independent. Then $e_\psi e_\vp=e_{\psi\vp}$ and $s^*_\vp
    s_\psi=s_\psi s^*_\vp$.
    \item [(b)] If $s^*_\vp s_\psi=s_\psi s^*_\vp$, then $\psi$ and
    $\vp$ commute and are independent.
\end{itemize}\elemma

\bproof (a) Since $e_\vp$ is exactly the orthogonal projection
  onto $\ell^2(\vp G)$, condition (c) in Lemma \ref{prime} translates to
  $e_\vp e_\psi =e_{\vp\psi}$ or $s_\vp s^*_\vp s_\psi s^*_\psi
  =s_\vp s_\psi s^*_\psi
  s^*_\vp = s_\psi s_\vp s^*_\psi
  s^*_\vp$ (for the last equality we use the fact that $\vp\psi =\psi\vp$).
  This is true if and only if $s^*_\vp s_\psi=
  s_\psi s^*_\vp$.

(b) Lemma \ref{vw} shows that the hypothesis implies that $s_\vp$
and $s_\psi$ and thus also $\vp$ and $\psi$, commute. Now we use
again that, for commuting $s_\vp,s_\psi$, the identity $s^*_\vp
s_\psi= s_\psi s^*_\vp$ holds if and only if $s_\vp s^*_\vp s_\psi
s^*_\psi =s_\vp s_\psi s^*_\psi s^*_\vp$, thus iff $e_\psi
e_\vp=e_{\psi\vp}$ which is equivalent to condition (c) in Lemma
\ref{prime}. \eproof

\section{Orbit partitions for endomorphisms and for pairs of independent
commuting endomorphisms}\label{orbit}

We briefly discuss here the orbit partition corresponding to an
endomorphism and to a pair of independent commuting endomorphisms of
a compact abelian group $H$. As in the previous sections we consider
only surjective endomorphisms $\alpha$ with finite kernel for which
the subgroup $L_{\alpha}=\bigcup_n \Ker {\alpha}^n$ is dense in $H$.

Let $\alpha$ and $\beta$ be two such endomorphisms which commute and
are independent in the sense of Definition \ref{inp}. From the
remark after Definition \ref{indt} we see that then every power of
$\alpha$ is independent from every power of $\beta$. Therefore, from
Lemma \ref{strong}, we obtain the following identities

\bgl\label{11}L_\alpha\cap L_\beta =\{1\}\quad L_\alpha
L_\beta=L_{\alpha\beta}\quad \alpha (L_\beta)=L_\beta\quad \beta
(L_\alpha)=L_\alpha\egl

\begin{definition}
Suppose that $S$ is a commutative semigroup of endomorphisms of the
compact abelian group $H$. The orbit of the point $x\in H$ with
respect to the semigroup $S$ is the set:
$$\orb_S(x)=\{y\in H:\exists \sigma_1,\sigma_2\in S,\; \sigma_1(x)=\sigma_2(y)\}$$
\end{definition}

In particular if $S=\{\alpha^n, n\in \Nz\}$ then we obtain the orbit
of the endomorphism $\alpha$
$$\orb_{\alpha}(x)=\{y: \exists n, m\in\Nz,\; {\alpha}^n(y)={\alpha}^m(x)\}$$

If $S$ is generated by two commuting endomorphisms $\alpha$ and
$\beta$ then the orbit of a point $x$ with respect to the semigroup
$S$ is
$$\orb_S(x)= \orb_{\alpha,\beta}(x) = \{y: \exists n, m,s,t \in\Bbb
N,\; {\alpha}^n\beta^s(y)=\alpha^m\beta^t(x)\}$$

It is clear that this definition defines a partition of the group
$H$ into orbits. This partition is called the orbit partition of the
semigroup $S$. Sometimes it is better to speak about the orbit
equivalence relation as a Borel subset of the product $H\times H$
(i.e. the set of all pairs of points which belong to the same
orbit).

For an exact endomorphism the orbit of each point is a dense
countable set; so this partition is not measurable in the usual
sense. This means that there is no measurable structure on the space
of orbits (no natural quotient under that partition). This fact is
usually expressed as the ergodicity of the orbit equivalence
relation.

The structure of an orbit for one exact endomorphism $\alpha$ is
easy to describe: One obviously has $\orb_\alpha(1)=L_{\alpha}$; for
a generic point $x\in H$, using surjectivity of $\alpha$ one can
choose a two sided sequence $z_k, k \in \Zz$ such that $z_0=x$, and
$z_{k+1}=\alpha(z_k)$ for all $k  \in \Zz$. Then the $\alpha$-orbit
is given by \bgl\label{12}\orb_{\alpha}(x)=\bigcup_{k\in \Zz}
z_kL_{\alpha}\egl

\begin{theorem}(\cite{Bow},\cite{Vershik})\label{BV} The orbit partition of an
arbitrary measure preserving  or nonsingular endomorphism of a
standard measure space is hyperfinite, i.e. a limit of a sequence of
monotonously increasing equivalence relations with finite classes.
\end{theorem}
This is a generalisation of Dye's theorem which says that the orbit
partition for a single nonsingular automorphisms is hyperfinite.
Connes-Feldman-Weiss \cite{CFW} proved that a nonsingular action of
any amenable group is hyperfinite. They also gave an alternative
proof of Theorem \ref{BV}, \cite[Corollary 13]{CFW}. Moreover,
Corollary 12 in \cite{CFW} also shows that the orbit partition for a
semigroup generated by finitely many commuting measure preserving
endomorphisms as in Theorem \ref{BV} is hyperfinite.\mn

In the case of a pair of independent commuting endomorphisms
hyperfiniteness can be derived from the hyperfiniteness of the orbit
partition for a single endomorphism very directly. From the
definition of the orbit partition for a commutative semigroup and
Lemma \ref{strong}, resp. formulas (\ref{11}) and (\ref{12}), we can
conclude that the orbit partition of the semigroup generated by two
commuting independent endomorphisms $\alpha$ and $\beta$ has the
form:

$$\orb_{\alpha,\beta}(x)=\bigcup_{y\in \orb_{\alpha}(x)}\orb_{\beta}(y)
\;=\bigcup_{y\in \orb_{\beta}(x)}\orb_{\alpha}(y)$$

This is a partition of product type. Since $L_\alpha\cap
L_\beta=\{1\}$, the intersection of an $\alpha$-orbit with a
$\beta$-orbit consists of at most one point.

\section{The C*-algebra of a rational polymorphism}

We consider now a pair $\vp,\psi$ of injective endomorphisms of the
(countable) abelian group $G$ and define the operators
$s_\vp,s_\psi$ on $\ell^2G$ by $s_\vp (\xi_g)=\xi_{\vp (g)}$ and
$s_\psi (\xi_g)=\xi_{\psi (g)}$. The fact that $s_\vp,s_\psi$ are
isometries satisfying $s_\vp s_\psi^*=s_\psi^*s_\vp$ encodes the
hypothesis that $\vp,\psi$ are injective, commute and are
independent, see Lemma \ref{pi}. This will be our assumption from
now on. Moreover we assume that $\vp$ and $\psi$ are exact, i.e.
that
$$\bigcap \vp^n (G) = \bigcap \psi^n (G) =\{0\}$$
and that $\vp$, $\psi$ have finite cokernel. We call this a standard
pair of endomorphisms from now on.

\bdefin\label{poly} We say that $\Pi$ is a rational polymorphism of
the compact abelian group $H$, if it is determined as in Definition
\ref{alge} by a correspondence $C_{\hat{\vp},\hat{\psi}}$, where
$\vp,\psi$ is a  standard pair of endomorphisms of the dual group
$G=\hat{H}$. We then denote $\Pi$ by $\vp /\psi$.\edefin

Each element of $G$ induces a unitary operator $u_g$ on $\ell^2G$
defined by $u_g(\xi_t)=\xi_{g+t}$. For a rational polymorphism as in
Definition \ref{poly}, $s_\vp$ commutes with $s^*_\psi$ and the
partial isometry $s=s_\varphi s^*_\psi$ describes the Markov
operator. We propose to describe the C*-algebra generated by $s$ and
the $u_g,\,g\in G$.

The operators $u_g:g\in G$ and $s=s_\vp s_\psi^*=s_\psi^* s_\vp$ in
$\cL (\ell^2G)$ satisfy the following relations:

\begin{itemize}
  \item[P1] $s^n$ is a partial isometry for each $n=1,2,\ldots$
  (i.e. the $s^{n*}s^n$ are projections).
  \item[P2] $g\mapsto u_g,\,g\in G$ is a representation of $G$ by
  unitaries.
  \item[P3] For each $g$, we have $s u_{\psi(g)} =u_{\vp (g)}s$.
  \item[P4] One has $$\sum_{g\in G/\vp (G)}u_gss^*u_{-g}=1\qquad
  \sum_{g\in G/\psi (G)}u_gs^*su_{-g}=1$$
\end{itemize}
The third condition implies that $u_{\vp(g)}$ commutes with $ss^*$
and $u_{\psi (g)}$ commutes with $s^*s$ for all $g\in G$. This
explains why, in the last condition, summation over $G/\vp (G)$
makes sense.

\bdefin Let $\Pi=\vp/\psi$ be a rational polymorphism as in
Definition \ref{poly}. We denote by $\Afp$ the universal C*-algebra
with generators $s$ and $u_g,\, g\in G$ satisfying the relations
P1,P2,P3,P4. \edefin

Much of the basic analysis of $\Afp$ can be developed in parallel
to the discussion in section \ref{end}.

\blemma (a) For each $n$, we have
$$\sum_{g\in G/\vp^n (G)}u_gs^ns^{*n}u_{-g}=1\qquad
  \sum_{g\in G/\psi^n (G)}u_gs^{*n}s^nu_{-g}=1$$
(b) Any two projections of the form $u_gs^ns^{*n}u_g^*$ and
$u_hs^{*m}s^m u_h^*$, $g,h\in G,\, n,m\in\Nz$, commute. \elemma

\bproof (a) We observe first that a product $ab$ of two partial
isometries is a partial isometry if and only if the range projection
of $b$ commutes with the support projection of $a$. Moreover
$s^ns^{*n}\leq s^ms^{*m}$ if $n\geq m$. Therefore we conclude from
property P1 that all projections of the form $s^ns^{*n}$ and
$s^{*m}s^m$ commute.

Since $\vp (G)+\psi (G)=G$ by \ref{prime} (a), and thus $G/\vp
G=\psi G/\vp G$, the identity $\sum_{g\in G/\vp
(G)}u_gss^*u_{g}^*=1$ from P4 can be rewritten as

$$\sum_{h\in G/\vp (G)}u_{\psi (h)}ss^*u_{\psi (h)}^*=1$$

Conjugating this by $s\,\cdot\,s^*$ and using property P3 one finds
that $ss^*=$

$\sum_{h\in G/\vp (G)}u_{\vp (h)}s^2s^{*2}u_{\vp (h)}^*$ and thus,
applying P4 again

$$\sum_{g,h\in G/\vp (G)}u_{g+\vp (h)}ss^*u_{g+\vp (h)}^*=
\sum_{g\in G/\vp^2 (G)}u_{g}s^2s^{*2}u_{g}^*=1$$

Iterating this procedure we find the first identity in (a). The
proof of the second identity is of course the same, replacing $\vp$
by $\psi$ and $s^ns^{*n}$ by $s^{*n}s^n$.

(b) Using the comment after Definition \ref{indt}, for any $n,m\in
\Nz$ we have $\vp^m(G)+\psi^n (G)=G$. Therefore we can rewrite any
two projections of the form $u_gs^{*n}s^nu_g^*$ and
$u_hs^ms^{*m}u_h^*$ as $u_{\vp^m(g')}s^{*n}s^nu_{\vp^m(g')}^*$ and
$u_{\psi^n(h')}s^ms^{*m}u_{\psi^n(h')}^*$. However, by P3,
$u_{\vp^m(g')}$ commutes with $s^ms^{*m}$ and $u_{\psi^n(h')}$
commutes with $s^{*n}s^n$. \eproof

By point (b) of this Lemma, the projections of the form
$u_gs^ns^{*n}u_{-g}$, $u_hs^{*m}s^{m}u_{-h}$ $g,h\in G,\,n,m\in\Nz$,
generate a commutative subalgebra $\cD$ of $\Afp$. The following
Lemma is then basically a special case of Lemma \ref{D}.

\blemma\label{DP} The C*-subalgebra $\cD$ of $\Afp$ generated by all
projections of the form $u_gs^ns^{*n}u_g^*,\,u_gs^{*m}s^m u_g^*$,
$g\in G, n,m\in \Nz$ is commutative. Its spectrum is the completion
$$G_{\psi\vp}=\mathop{\lim}\limits_{ {\scriptstyle\longleftarrow_n} }G/
(\psi\vp)^nG$$

$G$ acts on $\cD$ via $d\mapsto u_gdu_g^*$, $g\in G,\, d\in \cD$.
This action corresponds to the natural action of the dense subgroup
$G$ on its completion $G_{\psi\vp}$ via translation.

We have $G_{\psi\vp}\cong G_\vp\times G_\psi$ and $\cD\cong
\cD_\vp\otimes \cD_\psi$ with $\cD_\vp \cong C(G_\vp)$ and
$\cD_\psi\cong C(G_\psi)$.

The map $\cD s^*s\to \cD ss^*$ given by $x\mapsto sxs^*$ corresponds
to the map induced by $\vp\psi^{-1}$ on $(\psi (G))_{\psi\vp}\cong
G_\vp\times (\psi (G))_{\psi}$.\elemma

\bproof The first part of the assertion is just a special case of
Lemma \ref{D}. According to Remark \ref{chin} we have
$G_{\psi\vp}\cong G_\vp\times G_\psi$. \eproof

We will denote the compact abelian group $G_{\psi\vp}\cong
G_\psi\times G_\vp$ by $K$. Its dual group is the discrete abelian
group
$$L = \mathop{\lim}\limits_{ {\scriptstyle\longrightarrow_n}}\Ker
((\hat{\psi}\hat{\vp})^n :H\to H)$$ Because of the condition that we
impose on $\psi,\vp$, $L$ can be considered as a dense subgroup of
$H$.

The torus $\Tz$ acts on $\Afp$ by automorphisms $\alpha_t,\, t\in
\Rz$ defined by
$$\alpha_t (s)=e^{it}s\qquad \alpha_t (u_g)=u_g$$ The fixed point
algebra $B_{\psi\vp}$ is the subalgebra of $\Afp$ generated by all
$u_g,\,g\in G$ and by the $s^ns^{*n},\,s^{*n}s^n$. Integration over
$\Tz$ gives a faithful conditional expectation $\Afp\to B_{\psi\vp}$.

Applying Lemma \ref{B} to the endomorphism $\psi\vp$ we obtain
\blemma\label{Bfp} The C*-subalgebra $B_{\psi\vp}$ of $\Afp$ generated by
$C(H)$ together with $C(K)$ (or equivalently by $C^*G$ together with
$C^*L$) is isomorphic to the crossed product $C(K)\rtimes G$. It is
simple and has a unique trace. \elemma

As in section \ref{end} we also have the isomorphisms

$$B_{\psi\vp}\cong C(K)\rtimes G\cong C^*L\rtimes G\cong C^*G\rtimes L\cong
C(H)\rtimes L$$\mn Now, again as in section \ref{end}, we compose
the conditional expectation $\Afp\to B_{\psi\vp}$ with the natural
expectation $B_{\psi\vp}\to\cD$, coming from the crossed product
representation, to obtain a faithful conditional expectation
$E:\Afp\to\cD$.

\blemma\label{elem}\begin{itemize}
         \item[(a)] $s^*u_gs\neq 0\Rightarrow\; g\in \vp (G)$ and
         $su_gs^*\neq 0\Rightarrow\; g\in \psi (G)$
         \item[(b)] $su_gs\neq 0 \Rightarrow g\in \vp (G)\cap\psi
         (G)$
         \item[(c)] Every element in $\Afp$ can be
         approximated by linear combinations of elements of the
         form $u_gs^{*n}du_ts^mu_h$ or $u_gs^ndu_ts^{*m}u_h$ with
         $d\in \cD$, $n,m\geq 0$, $g,h,t\in G$.
       \end{itemize} \elemma

\bproof (a) follows immediately from condition P4. (b) follows from
the identities $$f_1=\sum_{g\in \vp (G)/(\psi (G)\cap\vp
(G))}u_gf_1f_2u_{-g}\qquad f_2=\sum_{h\in \psi (G)/(\psi (G)\cap\vp
(G))}u_hf_1f_2u_{-h}$$ where $f_1=ss^*,\,f_2=s^*s$.

Finally, using (a) and (b) one sees that the set of elements
described in (c) is invariant under multiplication by $s, s^*$,
$u_g$ on the right or on the left. \eproof

\blemma\begin{itemize}
         \item[(a)] $s$ and $s^*$ normalize $\cD$.
         \item[(b)] An element $z$ of the form $u_gs^{*n}du_ts^mu_h$ or
         $u_gs^ndu_ts^{*m}u_h$ with $d\neq 0$ is in $\cD$ if and only if
         $n=m$ and $g+h+t=0$.
         \item[(c)] Let $z$ be an element of the form $u_gs^{*n}du_ts^mu_h$
         or $u_gs^ndu_ts^{*m}u_h$. If $n\neq m$ or $g+h+t\neq 0$, then
         $E(z)=0$.
       \end{itemize}\elemma
\bproof (a) follows from Lemma \ref{elem} (a) and (b).

(b) If $n=m$ and $g+h+t=0$, then combining Lemma \ref{elem} (a) and
(b) with the fact that $s,s^*$ normalize $\cD$, we see that $z\in
\cD$. Conversely, if $z$ is in $\cD$, then $z$ has to fixed by
$\alpha_t,\,t\in\Rz$ whence $n=m$ and by $\beta_\chi,\,\chi\in
\hat{G}$ whence $g+h+t=0$. (c) is immediate from the definition of
$E$. \eproof

\blemma\label{cond} Let \bgl\label{z}z = d + \sum_{i=1}^k
u_{g_i}s^{*n_i}d_iu_{t_i}s^{m_i}u_{h_i}+ \sum_{i=k+1}^m
u_{g_i}s^{n_i}d_iu_{t_i}s^{*m_i}u_{h_i}\egl be an element of $\Afp$
(cf. Lemma \ref{elem} (c)) such that for each $i$, either $n_i\neq
m_i$ or $g_i+h_i+t_i\neq 0$ with $d,d_i\in \cD$. Then

\begin{itemize}
  \item[(a)] There is a projection $e\in
  \cD$, such that for all $i$ we have $e u_{g_i}s^{*n_i}d_iu_{t_i}s^{m_i}
  u_{h_i}e =0$ and $e u_{g_i}s^{n_i}d_iu_{t_i}s^{*m_i}u_{h_i}e =0$,
  and such
  that $e de=\lambda e$ with $\abs{\lambda}=\|E(z)\|$.
   \item[(b)] For the projection $e$ in (a) one
   has $eze=\lambda e$ with $\abs{\lambda}=\|E(z)\|$.
\end{itemize}\elemma

\bproof (a) $\cD$ is normalized by $s,s^*$ and the $u_g$. Arguing
exactly as in the proof of Theorem \ref{pure}, we see that there is
a projection $e$ in $\cD$ which is transported to an orthogonal
projection by each of the $u_{g_i}s^{*n_i}u_{t_i}s^{m_i}u_{h_i}$ and
$u_{g_i}s^{n_i}u_{t_i}s^{*m_i}u_{h_i}$ and such that $eze=\lambda e$
with $\abs{\lambda}=\|E(z)\|=\|d\|$. In fact for $e$ we can choose
$e'_1\otimes e'_2\in \cD_\vp\otimes \cD_\psi\cong \cD$, where $e'_1$
and $e'_2$ are chosen for $\vp$ and $\psi$ respectively as in the
proof of theorem \ref{pure}. One has $E(z)=d$, $eze=de$ and (b) is
then clear.\eproof

The algebra $B_{\psi\vp}$ is the inductive limit of the algebras
$B_n\cong M_{N(\psi\vp)^n}(C^*(G))$ (see the discussion at the
beginning of section \ref{Kend}). It naturally contains the inductive
limit $U_{\psi\vp}$ of the subalgebras $M_{N(\psi\vp)^n}(\Cz)
\subset M_{N(\psi\vp)^n}(C^*(G))$. This inductive limit is a $UHF$-algebra
of type $N(\psi\vp)^\infty$.

\bremark We have $N(\psi\vp)=N(\psi)N(\vp)$ (this is a slight generalization of the well known fact that the absolute norm in number theory is multiplicative). The proof follows from the facts that $G/\psi (G)\cong G/\psi\vp (G)\big{/}(\psi (G)/\psi\vp (G))$ and that $\psi (G)/\psi\vp (G)\cong G/\vp (G)$.\eremark

\blemma\label{t} Assume that $N(\vp)\geq N(\psi)$. There is a unitary $u\in M_{N(\psi\vp)}(\Cz)\subset U_{\psi\vp}$ such that the partial isometry $t=us$ satisfies $tt^*\leq t^*t$.

Let $\tau$ be the unique trace state on $B_{\psi\vp}$ and $x\in
B_{\psi\vp}$ such that $x=xt^*t$. Then $\tau
(txt^*)=N(\vp)^{-1}N(\psi)\tau (x)$. \elemma

\bproof The projections $ss^*$ and $s^*s$ are both contained in
$U_{\psi\vp}$ (in fact even in $M_{N(\psi\vp)}(\Cz)\subset
U_{\psi\vp}$) and satisfy $N(\vp)^{-1}=\tau (ss^*)\leq \tau
(s^*s)=N(\psi)^{-1}$ for the unique trace state $\tau$ on
$U_{\psi\vp}$. Therefore we can choose a unitary $u\in
M_{N(\psi\vp)}(\Cz)$ as required.

The unique trace state $\tau$ on $B_{\psi\vp}$ is obviously given as
the pointwise limit of the natural trace states $\tau_n$ on
$B_n\cong M_{N(\psi\vp)^n}(C^*(G))$ (obtained from the composition
of the natural map $M_{N(\psi\vp)^n}(C^*(G))\to
M_{N(\psi\vp)^n}(\Cz)$ given by the trivial representation of $G$
and the normalized trace on $M_{N(\psi\vp)^n}(\Cz)$). These states
clearly satisfy $\tau_n (sxs^*)=N(\vp)^{-1}N(\psi)\tau_n (xs^*s)$
for $x\in B_n$ (under the isomorphism $M_{N(\psi\vp)}(\Cz) \cong
M_{N(\psi)}(\Cz)\otimes M_{N(\vp)}(\Cz)$ we have $s(e_{11}\otimes
1)s^*=1\otimes e_{11}$). Therefore $\tau
(sxs^*)=N(\vp)^{-1}N(\psi)\tau (s^*sx)$.\eproof

\btheo\label{posi} Let $\Pi =\vp /\psi$ be a rational polymorphism
of the compact abelian group $H$. The associated universal
C*-algebra $\Afp$ is simple. It is purely infinite if $N(\vp)\neq
N(\psi)$. \etheo

\bproof Let $I$ be a non-zero closed ideal in $\Afp$ and $h$ a
non-zero positive element in $I$. Then $E(h)\neq 0$, and by Lemma
\ref{cond} in combination with lemma \ref{elem} (c), there is a
projection $e$ in $\cD$ (thus also in the subalgebra $B_{\psi\vp}$)
such that $ehe-\|E(h)\|e<\|E(h)\|/2$. It follows that $e\in I$.
Since $e\in B_{\psi\vp}$ and $B_{\psi\vp}$ is simple unital, also
$1\in I$. This shows that $\Afp$ is simple.

Possibly replacing $s$ by $s^*$ we can always assume that $N(\vp)\geq
N(\psi)$.

If now $N(\vp)> N(\psi)$, then choosing $t$ as in Lemma \ref{t}, we
have by \ref{t} that $\tau (t^kt^{*k})= (N(\vp)N(\psi)^{-1})^{-k}$.
Take $k$ large enough so that $\tau (t^kt^{*k})<\tau (e)$. In the
UHF-algebra $U_{\psi\vp}$ there is then a unitary $u$ such that
$ut^kt^{*k}u^*\leq e$. Thus $t^{*k}u^*hut^k$ is invertible in the
unital C*-algebra $t^*t\Afp t^*t$. This shows that $t^*t\Afp t^*t$
is purely infinite which immediately implies that $\Afp$ is purely
infinite too.\eproof

\bcor Let $\Pi =\vp /\psi$ be a rational polymorphism of the compact
abelian group $H$. The natural map from $\Afp$ to the C*-algebra of
operators on $L^2(H)$ generated by $s_\Pi$ and $C(H)$ is an
isomorphism.\ecor \bproof The map exists and is surjective by
universality of $\Afp$. It is injective since $\Afp$ is simple by
Theorem \ref{posi}.\eproof

\bcor\label{cross} Let $h=t^*t$ be the support projection of the
partial isometry $t$ defined in Lemma \ref{t} and $\gamma$ the
endomorphism of $hB_{\psi\vp}h$ defined by $\gamma (x)=txt^*$. Then
$h\Afp h$ is isomorphic to the semigroup crossed product
$B_{\psi\vp}\rtimes_\gamma \Nz$ (and thus $\Afp$ itself is Morita
equivalent to this crossed product).

If $N(\vp)=N(\psi)$, $\Afp$ has a unique trace state.\ecor

\bproof By definition, clearly there is a surjective map from the
universal algebra $h\Afp h$ onto the crossed product. By simplicity
of $\Afp$, this map is an isomorphism.

By Lemma \ref{t}, if $N(\psi)=N(\vp)$, then $\tau(txt^*)=t(x)$ for
$x$ in the hereditary subalgebra of $B_{\psi\vp}$ generated by
$t^*t$. Therefore $\tau$ extends uniquely to a trace state on
$\Afp$. This follows for instance from \cite{Nesh}, Corollary 1.2.
In fact, $\Afp$ can be seen as the C*-algebra of a groupoid with
object space $K$ and with the set of points in $K$ with non-trivial
stabilizer group of zero measure for $\tau$.

In our case at hand one can also give a very easy direct proof (cf.
the argument given in the proof of Lemma \ref{B}). In fact one can
easily construct partitions of unity in $C(K)$ consisting of
pairwise orthogonal projections $e_i$ with small trace such that for
an element $z$ such as in (\ref{z}) one has $e_ize_i=de_i$ for
nearly all $i$. We omit the details.\eproof

\section{Computation of the $K$-theory of $\Afp$}
From formula (\ref{KBf}) in section \ref{Kend} we conclude
\bgl\label{ind}K_*(B_{\psi\vp})=\mathop{\lim}\limits_{
\mathop{{\scriptstyle\longrightarrow}}
\limits_{b(\psi\vp)}}K_*(C^*(G))\egl

where $b(\psi\vp)$ is the homomorphism playing the role of the
$b(\vp)$ from section \ref{Kend} for the endomorphism $\psi\vp$
(i.e. $b(\psi\vp) = \kappa_1^{-1}\iota^{\psi\vp}\kappa_0$).

In the following we need however a somewhat finer analysis of the
$K$-theory homomorphisms involved in this formula. For this we
consider the C*-subalgebras $B_{mn}=C(G/\psi^m\vp^nG)\rtimes G\cong
M_{N(\psi^m)N(\vp^n)}(C^*G)$ of $B_{\psi\vp}$ and the corresponding
natural maps $\kappa_{mn}:K_*(C^*G) \to K_*(B_{mn})$.

As in section \ref{Kend} we denote by $\iota^\psi_{mn}$,
$\iota^\vp_{mn}$ the maps $K_*B_{mn}\to K_*B_{m+1,n}$ and
$K_*B_{mn}\to K_*B_{m,n+1}$ induced by the canonical inclusions and
set $b(\vp)_{mn}=\kappa_{m,n+1}^{-1}\iota^\vp_{mn}\kappa_{mn}$,
$b(\psi)_{mn}=\kappa_{m+1,n}^{-1}\iota^\psi_{mn}\kappa_{mn}$

\blemma\label{ab} The endomorphisms $b(\vp)_{mn}$ and $b(\psi)_{mn}$
of $K_*(C^*G)$ do not depend on $m,n$. Denoting these endomorphisms
by $b(\vp)$ and $b(\psi)$ we have $b(\psi\vp)
=b(\psi)b(\vp)=b(\vp)b(\psi)$.\elemma

\bproof The independence of $b(\vp)_{mn}$ and $b(\psi)_{mn}$ from
$m,n$ follows exactly as in section \ref{Kend}. The fact that
$b(\psi\vp)$ is the product of $b(\psi)$ and $b(\vp)$ (and that
these maps commute) can be seen from the following commutative
diagram

$$\xymatrix{K_*C^*(G)\ar[r]^{b(\vp)}\ar[d]^{\kappa_{00}}_=&K_*C^*(G)\ar[r]^
{b(\psi)}
\ar[d]^{\kappa_{01}}&K_*C^*(G)\ar[d]^{\kappa_{11}} \\
K_*B_{00}\ar[r]^{\iota^\vp_{00}}&K_*B_{01}\ar[r]^{\iota^\psi_{01}}&
K_*B_{11} }$$ and the fact that
$\iota^\psi_{01}\iota^\vp_{00}=\iota_0^{\psi\vp}$ (and the analogous
diagram with the order of $b(\vp)$ and $b(\psi)$ inverted, which
also commutes).\eproof

It will be convenient to represent $K_*(B_{\psi\vp})$ rather than as
in (\ref{ind}) as an inductive limit of the system $C_{mn}$ where
$C_{mn}=K_*(C^*G)$ for all $m,n\in\Nz$ and connecting maps
$b(\psi)^kb(\vp)^l:C_{mn}\to C_{(m+k)(n+l)}$.

The $(m,n)$ form a directed set for the order $(m,n)\leq
(m',n')\Leftrightarrow m\leq m'$ and $n\leq n'$. Lemma \ref{ab}
shows that $\mathop{\lim}\limits_{ \mathop{{\scriptstyle
\longrightarrow}}\limits_{b(\psi\vp)}}K_0(C^*(G))\cong
\mathop{\lim}\limits_{ {\scriptstyle\longrightarrow}_{mn}}C_{mn}$.
\mn

\btheo\label{Kfp} There is an exact sequence

$$\begin{array}{ccccc}
\mathop{\lim}\limits_{ \mathop{{\scriptstyle\longrightarrow}}
\limits_{b(\psi\vp)}}K_0(C^*(G)) &
\stackrel{\widetilde{b(\psi)}-\widetilde{b(\vp)}}{\lori} &
\mathop{\lim}\limits_{ \mathop{{\scriptstyle\longrightarrow}}
\limits_{b(\psi\vp)}}K_0(C^*(G)) &\lori & K_0\Afp
\\[3pt]
\uparrow  & & & & \downarrow
\\[2pt]
K_1\Afp & \longleftarrow & \mathop{\lim} \limits_{
\mathop{{\scriptstyle\longrightarrow}}
\limits_{b(\psi\vp)}}K_1(C^*(G))) & \stackrel{\widetilde{b(\psi)}
-\widetilde{b(\vp)}}{\longleftarrow} & \mathop{\lim} \limits_{
\mathop{{\scriptstyle\longrightarrow}}
\limits_{b(\psi\vp)}}K_1(C^*(G))
\end{array}$$\etheo

where $\widetilde{b(\psi)},\widetilde{b(\vp)}$ are the endomorphisms
of the inductive limit canonically induced by $b(\psi),b(\vp)$ (note
that $b(\psi),b(\vp)$ commute with $b(\psi\vp)$). \bproof From the
representation of $\Afp$ as a crossed product in Corollary
\ref{cross} we get as in section \ref{Kend} the Pimsner-Voiculescu
sequence \bgl\label{PV2}\xymatrix{K_*B_{\psi\vp}\ar[r]^{1-\gamma_*}&
K_*B_{\psi\vp}\ar[r]&K_*\Afp\ar@/^5mm/[ll]}\egl\mn As in section
\ref{Kend}, the maps $\kappa_{mn}:C_{mn}\to K_*(B_{mn})$ induce an
isomorphism $\kappa:\mathop{\lim}\limits_{
{\scriptstyle\longrightarrow}_{mn}}C_{mn} \lori K_*(B_{\psi\vp})$.
Let $\widetilde{b(\vp)}$ and $\widetilde{b(\psi)}$ denote the
automorphisms of $\mathop{\lim}\limits_{
{\scriptstyle\longrightarrow}_{mn}}C_{mn}$ induced by applying
$b(\vp)$ or $b(\psi)$, respectively to each $C_{mn}$ (thus
$\widetilde{b(\vp)}$ and $\widetilde{b(\psi)}$ correspond to the
left shifts with respect to $m$ or $n$, respectively, on a sequence
representing an element of the inductive limit). The identity
$\gamma_*\kappa_{mn}=\kappa_{(m+1)(n-1)}$ shows that $\gamma_*\kappa
= \kappa\widetilde{b(\psi)}^{-1}\widetilde{b(\vp)}$. Now
$\widetilde{b(\psi)}-\widetilde{b(\vp)}$ has the same kernel and
cokernel as $1-\gamma_*$. Therefore we get the asserted exact
sequence from (\ref{PV2}). \eproof

Because of exactness of the direct limit functor we obtain from Theorem
\ref{Kfp} the following short exact sequence

\bgl\label{sh}0\to\mathop{\lim}\limits_{
\mathop{{\scriptstyle\longrightarrow}} \limits_{b(\psi\vp)}}\;\Coker
(b(\psi) -b(\vp))\lori K_0\Afp\lori \mathop{\lim}\limits_{
\mathop{{\scriptstyle\longrightarrow}}
\limits_{b(\psi\vp)}}\Ker(b(\psi) -b(\vp))\to 0\egl

where the cokernel $(K_0(C^*(G))/(b(\psi) -b(\vp))K_0(C^*(G)))$ is
taken on $K_0$, while the kernel of $b(\psi) -b(\vp)$ is taken on
$K_1(C^*(G))$. Of course, there also is the analogous exact sequence
with the role of $K_0$ and $K_1$ interchanged.

\subsection{Example} Let $\vp$ and $\psi$ be endomorphisms of $\Zz$
determined by two relatively prime numbers $p$ and $q$. It is clear
that this pair of endomorphisms defines a rational polymorphism of
$\Tz$. The map $(\psi\vp)_*$ induced by the product is the identity
on $K_0(C(\Tz))=\Zz$ and multiplication by $pq$ on
$K_1(C(\Tz))=\Zz$. The identity $b(\psi\vp) (\psi\vp)_*=pq$ shows
that $b(\psi\vp)$ is multiplication by $pq$ on $K_0(C(\Tz))=\Zz$ and
the identity on $K_1(C(\Tz))=\Zz$. It follows that
$K_0(B_{\psi\vp})\cong \Zz[\frac {1}{pq}]$, $K_1(B_{\psi\vp})\cong
\Zz$.

Moreover $\Ker (b(\psi) -b(\vp))=0$ on $K_0(C(\Tz))$ and $\Ker
(b(\psi) -b(\vp))=\Zz$ on $K_1(C (\Tz))$.

$\Coker (b(\psi) -b(\vp))=\Zz/(p-q)$ on $K_0(C(\Tz))$ and $\Coker
(b(\psi) -b(\vp))=\Zz$ on $K_1(C (\Tz))$. From this and formula
(\ref{sh}), $K_*(\Afp)$ can easily be computed. Thus, if $p-q=1$, we
get $K_0\Afp=K_1\Afp=\Zz$. If for instance $p=5,\, q=3$, we get
$K_0\Afp=\Zz +\Zz /2,\;K_1\Afp=\Zz$.

\end{document}